\documentclass[12pt]{article}
\usepackage{latexsym,amsfonts,amssymb}
\newtheorem{thm}{Theorem}[section]

\newtheorem{defi}[thm]{Definition}
\newtheorem{rem}[thm]{Remark}
\newtheorem{exa}[thm]{Example}

\newtheorem{con}[thm]{Condition}

\date{}

\begin{document}
	\title{\bf  Random Almost Periodic Solutions of Random Dynamical
		Systems}
	\author{}
	\maketitle
	\centerline{Weili Zhang}\centerline{\small 1. Academy of Mathmatics and System Sciences, Chinese Academy of Sciences,
		Beijing, 100190, China}\centerline{\small
		2. School of Mathematical Sciences,
		University of Chinese Academy of Sciences,
		Beijing, 100049, China}	\centerline{\small
		E-mail: weilizhang17@mails.ucas.ac.cn}
	\vskip 1cm
	\centerline{Zuo-Huan Zheng} \centerline{\small 1. College of Mathmatics and Statistics,
		Hainan Normal University,}\centerline{\small
		\#99 Longkun Road, Haikou, 571158, China} \centerline{\small
		2. Institute of Applied Mathematics,
		Academy of Mathematics and System Sciences,}\centerline{\small
		Chinese Academy of Sciences,
		Beijing, 100190, China}\centerline{\small
		3. School of Mathematical Sciences,
		 University of Chinese Academy of Sciences,
		Beijing, 100049, China}
	\centerline{\small E-mail:
		zhzheng@amt.ac.cn}

	\vskip 0.5cm \noindent{\bf Abstract:}\quad In this paper, we first give the definition of random almost periodic solutions of random dynamical systems and give some examples. Then, we prove the existence of such random almost periodic solutions. Further, we introduce the definition of random almost periodic measure. Finally, we discuss the relationship between random almost periodic solutions and almost periodic measures.
	
	\smallskip
	
	\vskip 1cm
		
	
\section{Introduction}

\quad In the natural and social sciences, almost periodic phenomena are more easily seen than periodic phenomena. For example, many practical problems in celestial mechanics, power systems, ecological systems, economics, and engineering techniques can often be attributed to the search for periodic solutions and almost periodic solutions with ordinary differential equations as mathematical models, which for some problems(such as object motion and market supply and demand rule, etc.), investigating its almost periodic phenomena is sometimes more practical than examining its periodic phenomena. Therefore, it is of more practical significance to discuss its almost periodic properties(cf.\cite{MF,Y}).

\quad The theory of almost periodic function was founded by Bohr in the 1920s(cf.
\cite{B1, B2, B3}). Afterwards, a large amount of mathematical have made very important contributions in this respect, for example: V. Stepanov removes the requirement of continuity of functions and introduces Stapanov modules from Lebesgue function space, defining and studying a class of discontinuous almost periodic functions; Following this idea, H. Weyl and A. Besicoritch introduce a new semi-module to replace the classical uniform module, and obtain the extended definition of almost periodic function; Combining with the study of non-linear mechanics, N. Bogoliubov gives a profound analysis of Bohr theory from the idea of function approximation. Later developments are more closely related to ordinary differential equations, stability theory, dynamical systems, and partial differential equations(cf.\cite{B, AJ, ASB, CC}).

\quad In deterministic dynamical systems, almost periodic motions have been greatly developed due to the need to study practical problems and the development of other mathematical branches. Periodic motions are a special case of the almost periodic; almost periodic motions, in their turn, are a special case of recurrent motions. The mathematicians establishe the Lyapunov stability of almost periodic motion. They also investigate the question of when almost periodicity follows from Lyapunov stability(A.A.Markov has proved a stronger theorem)(cf.
\cite{NS}).

\quad In random dynamical systems, the concept of the stationary solution(random periodic solution) is a natural extension of the fixed point(periodic solution) in deterministic dynamical systems. In the past three decades, the stationary solution has developed vigorously(cf.\cite{A, AT, D, LZ}). In 2009, Zhao and zheng introduced the concept of random periodic solutions of random dynamical systems(cf.
\cite{ZZ}). Subsequently, a series of results related to stochastic periodic solutions and periodic probability solution are obtained(cf.
\cite{CFZ, CFLZ, FLZ, FWZ, FZZ, FZ, JQSY}). Recently, Sun and zheng introduce the weak random periodic solutions and weak random periodic measures of random dynamical system(cf.
\cite{SZ}). Some results of the almost periodic solutions(in distribution) to SDEs are obtained(cf.
\cite{BD, LLW, LW}). But some phenomena in life can not be accurately described by the almost periodic solution in distribution. In order to describe these phenomena better, we propose a random almost periodic solution of random dynamical systems.

\quad The paper is organized as follows. Section 2 is a preknowledge section. Section 3 gives definition of random almost periodic solution and gives some examples; Section 4 proves the sufficient conditions for the existence of random almost periodic solutions; Section 5 introduces the definition of almost periodic solution and discusses the relationship between random almost periodic solutions and almost periodic measures.

\section{preknowledge}	
\quad In this section, we will introduce the definitions of deterministic dynamical system, fixed point, periodic solution, almost periodic solution, random dynamical system, stationary solution and random periodic solution. Through this section, we assume that $(X,d)$ is a Polish space.

\quad A dynamical system or a flow on $X$ is a mapping $\psi:\mathbb{R}\times X \rightarrow X$ such that
\begin{equation}\label{1}
\psi(0,x)=x, \ \psi(t+s,x)=\psi(t, \psi(s,x)), \ \ for \ all \  x\in X \ and \ t\in \mathbb{R}.
\end{equation}

\quad For a deterministic dynamical system $\psi_t: X\rightarrow X$ over time $t\in \mathbb{R}$.\\
A fixed point is a point $x\in X$ such that
\begin{equation}\label{2}
\psi_t(x)=x, \ \ \ for \ all \  t\in \mathbb{R}.
\end{equation}
A periodic solution with periodic $\tau$ is a periodic function $\phi:I\rightarrow X$ such that
\begin{equation}\label{3}
\phi(t+\tau)=\phi(t), \psi(\phi(t_0))=\phi(t+t_0),\ \ \ for \ all \  t, t_0\in \mathbb{R}.
\end{equation}
And an almost periodic solution is an almost periodic function $\phi:\mathbb{R}\rightarrow X$ such that for any $\varepsilon>0$ there exists a number $l(\varepsilon)$ defining a relatively dense set of numbers $\{\tau_k\}$ which possess the following property:
\begin{equation}\label{4}
\rho(\phi(t+\tau_k),\phi(t))<\varepsilon, \psi(\phi(t_0))=\phi(t+t_0),\ \ \ for \ all \  t, t_0\in \mathbb{R}.
\end{equation}

\quad Let $(\Omega,{\cal F},P)$ be a probability space. A family of mapping(also called transformations) on the sample space, $\theta_t:\Omega\rightarrow\Omega, t\in \mathbb{R}$, is called a measurable dynamical system (or a measurable flow) if the following conditions are satisfied \\
\quad (i) Identity property: $\theta_0=Id$;\\
\quad (ii) Flow property: $\theta_{t+s}=\theta_{t}\circ \theta_{s}$;\\
\quad (iii) Measurability: $(\omega, t)\rightarrow \theta_{t}\omega$ is measurable.\\
It is called a measure-preserving dynamical system if, furthermore\\
\quad (iv) Measure-preserving proserving property: $P(\theta_{t}(A))=P(A)$, for every $A\in {\cal F}$ and $t\in \mathbb{R}$. \\
In this case, $P$ is called an invariant measure with respect to the dynamical system $\theta_t$.

\quad Consider an SDE system in $X$
\begin{equation}\label{5}
dX_t=b(X_t)dt+\sigma(X_t)dB_t
\end{equation}
In the definition of Brownian motion $B_t$, the probability space $(\Omega,{\cal F},P)$ is arbitrary. Now we introduce a specific or natural probability space in order to facilitate dynamical syatems research. We will call it canonical probability space for this SDE system. The canonical probability space for an SDE system (\ref{5}) in $X$ is $(\Omega,{\cal F},P)=(C(\mathbb{T},X),{\cal{B}}(C(\mathbb{T},X), P_B)$. The canonical sample space is now very concrete, as the samples are curves (i.e., paths of Brownian motion)(cf.\cite{KS}).

\quad The Wiener shift $\theta_t$ is defined as a mapping in the canonical sample space $\Omega$, for each fixed $t\in\mathbb{R}$,
\begin{equation}\label{6}
\theta_t:\Omega\rightarrow\Omega,\
\omega\mapsto\breve\omega \ \ \ such \  \ that \ \ \ \theta_{t}\omega (s)=\breve\omega(s):=\omega(t+s)-\omega(t), \ s,t\in\mathbb{R}.
\end{equation}

\quad By a simple calculation, we see that $\theta_0=Id$(the identity mapping in $\omega$) and $\theta_{s+t}=\theta_s\theta_t$. Moreover, $(w,t)\rightarrow \theta_{t}\omega$ is continuous and therefore measurable. Hence, the Wiener shift is a measurable dynamical system (or a flow) in$\Omega$. The above equation (\ref{6}) means that
\begin{equation}\label{7}
B_s(\theta_{t}\omega)=B_{t+s}(\omega)-B_t(\omega)
\end{equation}

\quad When $s$ is infinitesimally small, the right hand side is $dB_t(\omega)$. Thus $\theta_t$ is closely related to the noise in the stochastic system (\ref{4}) and is often called the driving flow.

\quad A measurable random dynamical system on the measurable space $(X,\cal B)$ over (or covering, or extending) a metric dynamical syatem $(\Omega,{\cal F},P,(\theta ({t})_{{t}\in\mathbb{R}}))$ with time $\mathbb{R}$ is a mapping $\Phi:{\mathbb{R}}\times\Omega\times X \rightarrow X, \  (t,\omega,x)\mapsto\Phi(t,\omega,x)$ with the following properties:\\
\quad (i) Measurability: $\Phi$ is $({\cal B}({\mathbb T}))\otimes {\cal F}\otimes {\cal B}(X), {\cal B}(X))$-measurable.\\	
\quad (ii) Cocycle property: The mappings $\Phi(t,\omega):=\Phi(t,\omega,\cdot) : X\rightarrow X$ form a cocycle over $\theta(\cdot)$, i.e. they satisfy
\begin{equation}\label{8}
\Phi(0,\omega)={id_X},  \ \ for \ all \  \omega\in \Omega \ \  (if \ 0\in \mathbb{R}),
\end{equation}
\begin{equation}\label{9}
\Phi(t+s,\omega)=\Phi(t,\theta_s\omega)\circ\Phi(s,\omega),  \ \ for \ all \  s,t\in\mathbb{R}, \ \omega\in \Omega.
\end{equation}	

\quad For random dynamical systems, it would be reasonable to say that stationary states are not actually steady states in the sence of (\ref{2}) and periodic solutions in the sense of (\ref{3}). Due to the fact that random exteral force pumps to the system constantly, the relation (\ref{2}) or (\ref{3}) breaks down.

\quad A random variable $y(\omega)$ is called a stationary solution (or random fixed point) for a random dynamical system $\Phi$ if
\begin {equation}\label{10}
\Phi(t,\omega,y(\omega))=y(\theta_t\omega) \ a.s., \ for \ all \ t\in\mathbb{R}.
\end{equation}

\quad A random periodic solution for a random dynamical system $\Phi$ is a $({\cal B}({\mathbb{R}}))\otimes{\cal B}(\Omega), {\cal B}(X))$-measurable map $\tilde{y}: \mathbb{R}\times\Omega\rightarrow X$ such that for almost all $\omega\in\Omega$,
\begin {equation}\label{11}
\Phi(t,\theta_s\omega)\tilde{y}(s,\omega)=\tilde{y}(t+s,\omega), \tilde{y}(s+\tau,\omega)=\tilde{y}(s,\theta_{\tau}\omega), \ \ for \ all \ \ t, s\in\mathbb{R}.
\end{equation}

\section{Random almost periodic solutions}

\quad Let $(\Omega, {\cal F}, P)$ be a probability space, $(X,d)$ be a compact metric space with metric $d$ and ${\cal B}(X)$ be it's Borel $\sigma$-algebra.
Consider a measurable random dynamical system $\Phi:\mathbb{R}\times\Omega\times X\rightarrow X$ on $(X, {\cal B}(X))$ over a metric dynamical system $(\Omega, {\cal F}, P, (\theta_t)_{t\in\mathbb{R}})$.

\quad Recall that a set $E\subset\mathbb{R}$ is said to be relatively dense if there exists a number $l>0$(length of indusion) such that every interval $[a,a+l]$, $a\in\mathbb{R}$, contains at least one point of $E$.

\begin{defi}
Given a measurable random dynamical system $\Phi:\mathbb{R}\times\Omega\times X\rightarrow X$ over $\theta$, an ${\cal F}$-measurable map $H:\mathbb {R}\times\Omega\rightarrow X$ is said to be a random almost periodic solution for the random dynamical system $\Phi$ if for any $\varepsilon>0$, there exists a number $l(\varepsilon)$ defining a relatively dense set of numbers $\{\tau_k\}$ shch that for almost all $\omega\in\Omega$,
	\begin{equation}\label{12}
	\Phi(t,\theta_s\omega)H(s,\omega)=H(t+s,\omega), \ d(H(s+\tau_k, \omega),H(s,\theta_{\tau_k}\omega))<\varepsilon, \ \ \forall t,s\in\mathbb{R}.
	\end{equation}	
\end{defi}

\begin{rem}
Stationary solutions and random periodic solutions are special cases of the random almost periodic solution.
\end{rem}

\begin{exa}
	Consider an SDE
	\begin{equation}\label{13}
	dX_t=-X_tdt+dB_t, \ \ \ \ X_0=x.
	\end{equation}
	This SDE defines a random dynamical system
	\begin{equation}\label{14}
	\Phi(t, \omega)x=\exp(-t)x+\int_{0}^{t} \exp(-(t-s))\, dB_s(\omega)
	\end{equation}	
	By examining the relation $\Phi(t, \omega)y(\omega)=y(\theta_t\omega)$ in this special case, a stationary solution of this random dynamical system is guessed out to be
	\begin{equation}\label{15}
	y(\omega)=\int_{-\infty}^{0}\exp(s)\, dB_s(\omega)
	\end{equation}
	Indeed, it follows from (\ref{14}) and (\ref{15}) that
	\begin{eqnarray}\label{16}
	\Phi(t, \omega)y(\omega)
	&=&\exp(-t)y(\omega)+\int_{0}^{t} \exp(-(t-s))\, dB_s(\omega)\nonumber\\
	&=&\exp(-t)\int_{-\infty}^{0}\exp(s)\, dB_s(\omega) + \int_{0}^{t} \exp(-(t-s))\, dB_s(\omega)\nonumber\\
	&=&\int_{-\infty}^{0}\exp(-(t-s))\, dB_s(\omega) + \int_{0}^{t} \exp(-(t-s))\, dB_s(\omega)\nonumber\\
	&=&\int_{-\infty}^{t} \exp(-(t-s))\, dB_s(\omega)
	\end{eqnarray}
	By (\ref{15}), we also see that
	\begin{eqnarray}\label{17}
	y(\theta_t\omega)&=&\int_{-\infty}^{0}\exp(s)\, dB_s(\theta_t\omega)\nonumber\\
	&=&\int_{-\infty}^{0}\exp(s)\, dB_{s+t}(\omega)\nonumber\\
	&=&\int_{-\infty}^{t} \exp(-(t-s))\, dB_s(\omega)
	\end{eqnarray}
	Thus, $\Phi(t, \omega)y(\omega)=y(\theta_t\omega)$, i.e. $y(\omega)=\int_{-\infty}^{0}\exp(s)\, dB_s(\omega)$ is a stationary solution for the random dynamical system (\ref{14}). \\
	Set $H(t, \omega)=y(\theta_t\omega)$, then we have
	\begin{eqnarray}\label{18}
	\Phi(t, \theta_s\omega)H(s, \omega)
	&=&\Phi(t, \theta_s\omega)y(\theta_s\omega)\nonumber\\
	&=&y(\theta_t\theta_s\omega)\nonumber\\
	&=&y(\theta_{t+s}\omega)\nonumber\\
	&=& H(t+s, \omega)
	\end{eqnarray}
	On the other hand, it is clearly that
	\begin{equation}\label{19}
	H(t+\tau, \omega)=H(t, \theta_{\tau}\omega), \ \ \forall \tau\in\mathbb{R}.
	\end{equation}
	So by (\ref{12}), $H(t, \omega)=y(\theta_t\omega)$ is a random almost periodic solution for the random dynamical system generated by SDE (\ref{14}).
\end{exa}

\begin{exa}
	Consider the stochastic equation
	\begin{equation}\label{20}
	dX_t=(\frac{3}{2}X_t-{X_t}^3)dt+X_tdB_t, \ \ \ \ X_0=x.
	\end{equation}
	The random dynamical system generated by this SDE is
	\begin{equation}\label{21}
	\Phi(t, \omega)x=\frac{x\exp({t+B_t(\omega)})}{{(1+2x^2\int_{0}^{t} \exp(2s+2B_s(\omega))\, ds)}^\frac{1}{2}}
	\end{equation}
	By examining the relation $\Phi(t, \omega)y(\omega)=y(\theta_t\omega)$ in this special case, it appears that a stationary solution of this random dynamical system is
	\begin{equation}\label{22}
	y(\omega)={(\int_{-\infty}^{0}\exp(2s+2B_s(\omega))\, ds)^{-\frac{1}{2}}}
	\end{equation}
	We can verify that this is indeed a stationary solution for this random dynamical system. Similarly to Example 3.3, set $H(t, \omega)=y(\theta_t\omega)$, we have
	\begin{equation}\label{23}
	\Phi(t,\theta_s\omega)H(s,\omega)=H(t+s,\omega), \
	H(t+\tau, \omega)=H(t, \theta_{\tau}\omega). \ \ \forall \tau\in\mathbb{R}
	\end{equation}
	So by (\ref{12}), $H(t, \omega)=y(\theta_t\omega)$ is a random almost periodic solution for the random dynamical system generated by SDE (\ref{21}).
\end{exa}	

\begin{exa}
	Consider the example given by Zhao and Zheng (see \cite[Section 2]{ZZ}).
	\begin{eqnarray}\label{24}
	\left\{
	\begin{array}{l}
	dx(t)=(x-y-x(x^2+y^2))dt+x\circ dB_t,\\
	dy(t)=(x+y-y(x^2+y^2))dt+y\circ dB_t.
	\end{array}
	\right.
	\end{eqnarray}	
	Here $B_t$, $t\in\mathbb{R}$ is a one-dimensional two-sided Brownian motion on the path space\\ $(C((-\infty,\infty)), {\cal B}(C((-\infty,\infty))),P)$
	with P-preserving map $\theta$ being  taken to the shift operator $(\theta_t\omega)(s)= \omega(t+s)-\omega(t)$ for $s,t\in \mathbb{R}$. Using polar	coordinates
	$$
	x=\rho\cos(2\pi\alpha),\ \ y=\rho\sin(2\pi\alpha),
	$$
	then we transform equation (\ref{24}) on $\mathbb{R}^2$ to the following equation on the cylinder $S^1\times\mathbb{R}$:
	\begin{eqnarray}\label{25}
	\left\{
	\begin{array}{l}
	d\rho(t)=(\rho(t)-\rho^3(t))dt+\rho(t)\circ dB_t,\\
	d\alpha(t)=\frac{1}{2\pi}dt.
	\end{array}
	\right.
	\end{eqnarray}	
	This SDE defines a random dynamical system
	\begin{equation}\label{26}	
	\Phi(t,\omega)(\alpha_0,\rho_0)=(\{\alpha_0+\frac{t}{2\pi}\}, \frac{\rho_0\exp({t+B_t(\omega)})}{(1+2\rho^2_0\int_{0}^{t} \exp(2s+2B_s(\omega))\, ds)^\frac{1}{2}})
	\end{equation}
	where the symbol $\{x\}$ is defined by $-\frac{1}{2}<\{x\}\leq \frac{1}{2}$ and $x=\{x\}+k$, $k$ is integer. By examing the relation $\Phi(t,\theta(s)\omega)Y(s,\omega)=Y(t+s,\omega)$, $Y(s+\tau,\omega)=Y(s,\theta_{\tau}\omega)$ in this special case, it appears that a random periodic solution for this random dynamical system is
	\begin{equation}\label{27}
	Y(t,\omega)=( \{\alpha_0+\frac{t}{2\pi}\},  {(\int_{-\infty}^{0}\exp(2s+2B_s(\theta_t\omega))\, ds)^{-\frac{1}{2}}})
	\end{equation}
	Indeed, it follows from (\ref{26}) and (\ref{27})
	\begin{eqnarray}\label{28}
	\Phi(t,\theta(s)\omega)Y(s,\omega)
	&=&(\{\alpha_0+\frac{s+t}{2\pi}\},\frac{Y(s,\omega)\exp({t+B_t(\theta(s)\omega)})}{{(1+2\rho^2_0\int_{0}^{t} \exp(2s+2B_s(\theta(s)\omega))\, ds)}^\frac{1}{2}})\nonumber\\
	&=&(\{\alpha_0+\frac{t+s}{2\pi}\}, {(\int_{-\infty}^{0}\exp(2s+2B_s(\theta_{t+s}\omega))\, ds)^{-\frac{1}{2}}})\nonumber\\
	&=&Y(t+s,\omega)
	\end{eqnarray}
	By (\ref{27}), we also see that
	\begin{eqnarray}
	Y(2\pi+t,\omega)
	&=&(\{\alpha_0+\frac{t+2\pi}{2\pi}\}, {(\int_{-\infty}^{0}\exp(2s+2B_s(\theta_{t+2\pi}\omega))\, ds)^{-\frac{1}{2}}})\nonumber\\
	&=&Y(t,\theta_{2\pi}\omega)
	\end{eqnarray}\label{29}
	Thus $Y(t,\omega)=( \{\alpha_0+\frac{t}{2\pi}\},  {(\int_{-\infty}^{0}\exp(2s+2B_s(\theta_t\omega))\, ds)^{-\frac{1}{2}}})$ is a random almost periodic solution for the random dynamical system. Set $H(t,\omega)=Y(t,\omega)$, then we have $$\Phi(t,\theta_s\omega)H(s,\omega)=H(t+s,\omega), \ H(s+2k\pi,\omega)=H(s,\theta_{2k\pi}\omega),\ k=0,\pm 1,\pm 2,\pm 3,\cdots$$
	So $H(t,\omega)$ is a random almost periodic solution with relatively dense set $\{2k\pi,k=0,\pm 1,\pm 2,\pm 3,\cdots\}$ for the random dynamical system generated by SDE (\ref{24}).
\end{exa}

\section{Existence of random almost periodic solutions}
\quad In this section, we will give a sufficient theorem for the existence of random almost periodic solutions. Consider a continuious time differentiable random dynamical syatem $\Phi$ over a metric dynamical system $(\Omega, {\cal F}, P, \theta)$ on $S^1\times S^1 \times \mathbb{R}^d$.

\begin{defi}
A $C^1$ perfect cocycle is a $({\cal B}(\mathbb{R})\otimes {\cal F}\otimes {\cal B}(S^1\times S^1\times \mathbb{R}^d), {\cal B}(S^1\times S^1\times \mathbb{R}^d))$-measurable random field $\Phi:\mathbb{R}\times \Omega\times S^1\times S^1\times \mathbb{R}^d\rightarrow S^1\times S^1\times \mathbb{R}^d$ satisfying the following conditions:\\
	(i) for each $\omega\in\Omega$, $\Phi(0,\omega)=Id$;\\
	(ii) for each $\omega\in\Omega$, $\Phi(t_1+t_2, \omega)x=\Phi(t_2, \theta_{t_1}\omega)\Phi(t_1, \omega)x$ for all $x\in S^1\times S^1\times \mathbb{R}^d$ and $t_1, t_2\in \mathbb{R}$;\\
	(iii) for each $\omega\in\Omega$, the mapping $\Phi(\cdot, \omega)\cdot: \mathbb{R}\times S^1\times S^1\times \mathbb{R}^d\rightarrow S^1\times S^1\times \mathbb{R}^d$ is continuous;\\
	(iv) for each $(t, \omega)\in\mathbb{R}\times \Omega$, the mapping $\Phi(t, \omega)\cdot: S^1\times S^1\times \mathbb{R}^d\rightarrow S^1\times S^1\times \mathbb{R}^d$ is a $C^1$ diffeomorphism.
\end{defi}

\begin{con}
	Assume that there exists two rational independent numbers $t_1$ and $t_2$, such that for any $(x, y, z)\in S^1\times S^1\times \mathbb{R}^d$,
	$$\Phi(t_1, \omega)(x, y, z)=(x, y^\prime, z^\prime)$$
	and
	$$\Phi(t_2, \omega)(x, y, z)=(x^\prime, y , z^\prime)$$
\end{con}

Note that $t_1$ is the time that a particle on $S^1$ rotate a full circle in $x$-direction and $t_2$ is the time that a particle on $S^1$ rotate a full circle in $y$-direction. Consider the following random system
\begin{eqnarray}\label{30}
\left\{
\begin{array}{l}
(x,y)=f(t;(x_0,y_0))\\
z=g(t, \omega, z_0).
\end{array}
\right.
\end{eqnarray}
where $(x_0, y_0)\in S^1\times S^1$ and $z_0\in \mathbb{R}^d$.
Since $t_1$ and $t_2$ are two rational independent numbers, for any $\varepsilon>0$ we can find a relatively dense set $\{\tau_k\}$ such that\\
$$d(f(\tau_k;(x,y)) \ mod \ \ \mathbb{Z}\times\mathbb{Z}, (x,y))<\frac{1}{\sqrt{2}}\varepsilon.$$
\begin{thm}
	Assume Condition 4.2 holds and $d(g(t+\tau_k, \omega, z), g(t, \theta_{\tau_k}\omega, z))<\frac{1}{\sqrt{2}}\varepsilon$ for all $z\in\mathbb{R}^d$, the perfect cocycle $\Phi:\mathbb{R}\times \Omega\times S^1\times S^1\times \mathbb{R}^d\rightarrow S^1\times S^1\times \mathbb{R}^d$ generated by (\ref{30}) has a random almost periodic solution.
\end{thm}
\noindent {\bf Proof.}\ \
Denote $H(s, \omega)=\Phi(s,\omega)(x_0,y_0,z_0)$, by Condition 4.2(ii), it's easy to see that $$\Phi(t,\theta_s\omega)H(s,\omega)=H(t+s,\omega).$$
By assumption, we also see that
\begin{eqnarray}\label{31}
&&d(H(t+\tau_k,\omega),H(t,\theta_{\tau_k}\omega))\nonumber\\
&=&\sqrt{d^2(f(\tau_k;(x_0,y_0)) \ mod \ \ \mathbb{Z}\times\mathbb{Z},,(x_0,y_0))+d^2(g(t+\tau_k, \omega, z_0), g(t, \theta_{\tau_k}\omega, z_0))}\nonumber\\
&<&\varepsilon.
\end{eqnarray}
Thus $H(s,\omega)$ is a random almost periodic solution for the random dynamical system $\Phi$.
\hfill\fbox\\

\quad  Next, we give an example of random almost periodic solution. In deterministic cases, consider an example of almost periodic motion on the torus $T^2(r,\alpha)$: $0\leq r<1$, $0\leq\alpha<1$. $(r+k,\alpha+k')\equiv(r,\alpha)$ if $k$ and $k'$ are integers.
\begin{eqnarray}\label{32}
\left\{
\begin{array}{l}
dr(t)=dt\\
d\alpha(t)=\gamma dt.
\end{array}
\right.
\end{eqnarray}
where $\gamma$ is irrational, in addition we define the distance between the points $(r_1,\alpha_1)$ and $(r_2,\alpha_2)$ as
\begin{equation}\label{33}
d((r_1,\alpha_1),(r_2,\alpha_2))=[\{r_1-r_2\}^2+\{\alpha_1-\alpha_2\}^2]^{\frac{1}{2}}
\end{equation}
(the symbol $\{x\}$ is defined by $-\frac{1}{2}<\{x\}\leq \frac{1}{2}$ and $x=\{x\}+k$, $k$ is integer).
By the definition of the almost periodic solution, this almost periodic solution genereted by (\ref{32}) implies the existence of a relatively dense set of numbers $\{\tau_k\}$, satisfying the condition: for any $\varepsilon>0$
\begin{equation}\label{34}
d((r(t+\tau_k),\alpha(t+\tau_k)),(r(t),\alpha(t)))<\varepsilon \ \ \forall t\in\mathbb{R}.
\end{equation}
i.e. $[\{r_1-r_2\}^2+\{\alpha_1-\alpha_2\}^2]^{\frac{1}{2}}<\varepsilon$.\\

\begin{exa}
	Consider the SDE on $\mathbb{R}^2\times S^1$
	\begin{eqnarray}\label{35}
	\left\{
	\begin{array}{l}
	dx_t=(\frac{3}{2}x_t-y_t-(x^2_t+y^2_t)x_t)dt+2\pi x_tdB_t\\
	dy_t=(x_t+\frac{3}{2}y_t-(x^2_t+y^2_t)y_t)dt+2\pi y_tdB_t\\
	dz_t=\gamma dt.
	\end{array}
	\right.
	\end{eqnarray}
	where $\gamma$ is irrational, $B_t$ is a one-dimensional standard Brownian motion. By the coordinate transformation $$x=rcos2\pi\alpha, \ y=x=rsin2\pi\alpha, \ z=z$$	
	we can transform (\ref{35}) on $\mathbb{R}^2\times S^1$ to the following SDE on 	$\mathbb{R} \times S^1\times S^1$
	\begin{eqnarray}\label{36}
	\left\{
	\begin{array}{l}
	dr_t=(\frac{3}{2}r_t-r^3_t)dt+r_tdB_t\\
	d\alpha_t=dt\\
	dz_t=\gamma dt.
	\end{array}
	\right.
	\end{eqnarray}
	The random dynamical system generated by this SDE is
	\begin{equation}\label{37}
	\Phi(t, \omega)(r_0,\alpha_0, z_0)=(\frac{r_0\exp({t+B_t(\omega)})}{(1+2r^2_0\int_{0}^{t} \exp(2s+2B_s(\omega))\, ds)^\frac{1}{2}},\{\alpha_0+t\}, \{z_0+\gamma t\})
	\end{equation}
	By examining the relation: for any $\varepsilon>0$, there exists a number $l(\varepsilon)$ defining a relatively dense set of numbers $\{\tau_k\}$ shch that for almost all $\omega\in\Omega$, $$\Phi(t, \theta(u)\omega)H(u,\omega))=H(u+t,\omega), \
	d(H(t+\tau_k,\omega),H(t,\theta(\tau_k)\omega))<\varepsilon \ \ \forall t,u\in\mathbb{R}$$ in this special case, a random almost periodic solution of this random dynamical system is gaussed out to be
	\begin{equation}\label{38}
	H(u,\omega)=((2\int_{-\infty}^{0}\exp(2s+2B_s(\theta_u\omega))\, ds)^{-\frac{1}{2}},\{\alpha_0+u\}, \{z_0+\gamma u\})
	\end{equation}
	Indeed, it follows from (\ref{37}) and (\ref{38})
	\begin{equation}\label{39}
	\Phi(t, \theta_u\omega)H(u,\omega)
	=(\frac{r_u\exp({t+B_t(\theta_u\omega)})}{(1+2r^2_u\int_{0}^{t} \exp(2s+2B_s(\theta_u\omega))\, ds)^\frac{1}{2}},\{\alpha_0+u+t\},\{z_0+\gamma(u+t)\})
	\end{equation}
	where $r_u=(2\int_{-\infty}^{0}\exp(2s+2B_s(\theta_u\omega))\, ds)^{-\frac{1}{2}}$,
	it's easy to see that
	\begin{eqnarray}\label{40}
	&&\frac{r_u\exp({t+B_t(\theta_u\omega)})}{(1+2r^2_u\int_{0}^{t} \exp(2s+2B_s(\theta_u\omega))\, ds)^\frac{1}{2}}\nonumber\\
	&=&\frac{\exp({t+B_t(\theta_u\omega)})}{(r^{-2}_u+2\int_{0}^{t} \exp(2s+2B_s(\theta_u\omega))\, ds)^\frac{1}{2}}\nonumber\\\nonumber\\
	&=&\frac{\exp({t+B_t(\theta_u\omega)})}{(2\int_{-\infty}^{0}\exp(2s+2B_s(\theta_u\omega))\, ds+2\int_{0}^{t} \exp(2s+2B_s(\theta_u\omega))\, ds)^\frac{1}{2}}\nonumber\\
	&=&\frac{\exp({t+B_t(\theta_u\omega)})}{(2\int_{-\infty}^{t} \exp(2s+2B_s(\theta_u\omega))\, ds)^\frac{1}{2}}\nonumber\\
	&=&(2\int_{-\infty}^{0}\exp(2s-2t+2B_s(\theta_u\omega)-2B_t(\theta_u\omega))\, ds)^{-\frac{1}{2}}\nonumber\\
	&=&(2\int_{-\infty}^{0}\exp(2s+2B_s(\theta_{u+t}\omega))\, ds)^{-\frac{1}{2}}
	\end{eqnarray}
	Thus,
	\begin{eqnarray}\label{41}
	\Phi(t, \theta_u\omega)H(u,\omega)
	&=&((2\int_{-\infty}^{0}\exp(2s+2B_s(\theta_{u+t}\omega))\, ds)^{-\frac{1}{2}},\{\alpha_0+u+t\}, \{z_0+\gamma(u+t)\})\nonumber\\
	&=&H(u+t,\omega)
	\end{eqnarray}
	By (\ref{37}), we also see that
	\begin{equation}\label{42}
	d(H(t+\tau,\omega),H(t,\theta(\tau)\omega))=[\{\tau\}^2+\{\gamma\tau\}^2]^{\frac{1}{2}}
	\end{equation}
	So for any $\varepsilon>0$, by (\ref{24}) there exists a relatively dense set $\{\tau_k\}$(same with the $\{\tau_k\}$ in (\ref{34})), s.t. $d(H(t+\tau_k,\omega),H(t,\theta(\tau_k)\omega))<\varepsilon, \ \ \forall t\in\mathbb{R}, \tau_k\in\{\tau_k\}$. Thus $H(u,\omega)$ defined by (\ref{38}) is a random almost periodic solution for the random dynamical system (\ref{37}).
\end{exa}

\quad Here we give a more general existence theorem.

\begin{thm}
	Let $X$ be a Polish space with metric $d$ and $\Phi $ be a measurable random dynamical system satisfying
	\begin{equation}\label{43}
	d(\Phi(t,\omega)x,\Phi(t,\omega)y)\leq Cd(x,y), \ \ \forall t\in\mathbb{R}
	\end{equation}
	for some constant C and $H:\mathbb {R}\times\Omega\rightarrow X$ be an ${\cal F}$-measurable map such that for any $\varepsilon>0$, there exists a number $l(\varepsilon)$ defining a relatively dense set $\{\tau_k\}$ such that for almost all $\omega\in\Omega$
	\begin{equation}\label{44}
	d(Y(0,\omega),\Phi(\tau_k,\theta_{-\tau_k}\omega)H(0,\theta_{-\tau_k}\omega))<\varepsilon
	\end{equation}
	Then $H(t,\omega)$ given by
	\begin{equation}\label{45}
	H(t,\omega):=\Phi(t,\omega)H(0,\omega)
	\end{equation}
	is a random almost periodic solution for the random dynamical system $\Phi$.
\end{thm}
\noindent {\bf Proof.}\ \
By the definition of the random dynamical system $\Phi$ and ${\cal F}$-measurable map $H$. For any $t,s\in\mathbb{R}$, almost all $\omega\in\Omega$,
\begin{eqnarray}\label{46}
\Phi(t,\theta_s\omega)H(s,\omega)
&=&\Phi(t,\theta_s\omega)\Phi(s,\omega)H(0,\omega)\nonumber\\
&=&\Phi(t+s,\omega)H(0,\omega)\nonumber\\
&=&H(t+s,\omega)
\end{eqnarray}
By (\ref{45}), $\Phi(\tau_k,\theta_{-\tau_k}\omega)H(0,\theta_{-\tau}\omega)=H(\tau_k,\theta_{-\tau_k})$. So for any $s\in\mathbb{R}$, by (\ref{46}) and (\ref{43})
\begin{eqnarray}\label{47}
d(H(s+\tau_k,\theta_{-\tau_k}\omega), H(s,\omega))
&=&d(\Phi(s,\omega)H(\tau_k,\theta_{-\tau_k}\omega),\Phi(s,\omega)H(0,\omega))\nonumber\\
&\leq& Cd(H(\tau_k,\theta_{-\tau_k}\omega,H(0,\omega))
\end{eqnarray}
For any $\varepsilon_1>0$, by (\ref{44}) there exists a number $l(\varepsilon_1)$ defining a relatively dense set $\{\tau_k\}$ such that for almost all $\omega\in\Omega$, $$d(H(s+\tau_k,\theta_{-\tau_k}\omega), H(s,\omega))<C\varepsilon_1$$
Therefore, $H$ is a random almost periodic solution for the random dynamical system $\Phi$.
\hfill\fbox\\

\section{Almost periodic probability measure}
Let $\Phi:\mathbb{R}\times\Omega\times X\rightarrow X$ be a measurable random dynamical system over $\theta$, consider a standard product measurable space $(\overline\Omega, \overline {\cal F})=(\Omega\times X, {\cal F}\otimes{\cal B}(X))$ and the skew-product of the metric dynamical system $(\Omega,{\cal F}, P, (\theta(t))_{t\in\mathbb{R}})$ and the cocycle $\Phi(t,\omega)$ on $X$, ${\Theta}(t):\overline\Omega\rightarrow\overline\Omega$,
\begin{equation}\label{48}
{\Theta}(t)(w,x)=(\theta_t\omega,\Phi(t,\omega)x), \ \ t\in\mathbb{R}.
\end{equation}

Denote
$$C_{BL}(\Omega\times X):=\{f\in C(\Omega\times X): ||f||=\sup\limits_{(\omega,x)\in\Omega\times X}|f(\omega,x)|; \  ||f(\omega,x)-f(\omega,y)||\leq C_1d(x,y), \forall \omega\in\Omega \}.$$

$$C_{BL}(X):=\{g\in C(X): ||g||=\sup\limits_{x\in X}|g(x)|; \  ||g(x)-g(y)||\leq C_2d(x,y)\}.$$
where $C_1$, $C_2$ are constants in $\mathbb{R}^+$.

Recall \\
\centerline {${\cal P}_P(\Omega\times X):=\{\mu:$ probability\ measure\ on $(\Omega\times X, {\cal F}\otimes{\cal B}(X))$ with marginal $P$ on $(\Omega, {\cal F})\}$}

and

\centerline {${\cal P}(X):=\{\lambda:$ probability\ measure\ on $(X, {\cal B}(X))\}$}

We endow ${\cal P}_P(\Omega\times X)$ with the $\rho$ metric:
\begin{equation}\label{49}
\rho(\mu,\nu)=sup\{|\mu(f)-\nu(f)|:f\in C_{BL}(\Omega\times X),||f||\leq 1 \}, \ \forall \mu,\nu\in {\cal P}_P(\Omega\times X).
\end{equation}

and ${\cal P}(X)$ with the $\rho_1$ metric:
\begin{equation}\label{50}
\rho_1(\mu,\nu)=sup\{|\mu(g)-\nu(g)|:g\in C_{BL}(X),||g||\leq 1 \}, \ \forall \mu,\nu\in {\cal P}(X).
\end{equation}

Suppose $\mu\in {\cal P}_P(\Omega\times X)$. We call a function $\mu_{\cdot}(\cdot):\Omega\times{\cal B}(X)\rightarrow [0,1]$ a factorization of $\mu$ with respect to $P$ if: \\
(i)for all $B\in{\cal B}(X)$, $\omega\mapsto\mu_{\omega}(B)$ is ${\cal F}$-measurable;\\
(ii)for $P-a.a.\omega\in\Omega$, $B\mapsto\mu_{\omega}(B)$ is a probability measure on $(X,{\cal B}(X))$;\\
(iii)for all $A\in{\cal F}\otimes{\cal B}(X)$
\begin{equation}\label{51}
\mu(A)=\int_{\Omega}\int_X \chi_A(\omega,x)\mu_{\omega}(dx)P(d\omega).
\end{equation}
Introducing the section $A_{\omega}:=\{x:(\omega,x)\in A\}$, (\ref{51}) can be written as
$$\mu(A)=\int_{\Omega}\mu_{\omega}(A_{\omega})P(d\omega).$$

\begin{defi}
	Given a measurable random dynamical system $\Phi:\mathbb{R}\times\Omega\times X\rightarrow X$ over $\theta$, a measure $\mu:\mathbb{R}\rightarrow{\cal P}_P(\Omega\times X)$ is called an almost periodic probability measure on $(\Omega\times X, {\cal F}\otimes{\cal B}(X))$ if for any $\varepsilon>0$, there exists a number $l(\varepsilon)$ defining a relatively dense set of numbers $\{\tau_k\}$ such that for almost all $\omega\in\Omega$,
	\begin{equation}\label{52}
	\Theta(t)\mu_s=\mu_{t+s}, \ \rho(\mu_{\tau_k+s},\mu_s)<\varepsilon, \ \forall t, s\in\mathbb{R}.
	\end{equation}
\end{defi}

\begin{thm}
	Let $\Phi:\mathbb{R}\times\Omega\times X\rightarrow X$ be a measurable random dynamical system over $\theta$, if it has a random  almost periodic solution $H:\mathbb{R}\times\Omega\rightarrow X$ with relatively dense set $\{\tau_k\}$. Then it has an almost periodic measure $\mu:\mathbb{R}\rightarrow{\cal P}_P(\Omega\times X)$ on $(\Omega\times X, {\cal F}\otimes{\cal B}(X))$. For any $f\in C_{BL}(\Omega\times X)$,
	\begin{equation}\label{53}
	\mu_t(f)=\int_{\Omega}f(\theta_{t}\omega,H(t,\omega))P(d\omega)
	\end{equation}	
	and for any $A\in{\cal F}\times{\cal B}(X)$,	
	\begin{equation}\label{54}
	\mu_t(A)=\int_{\Omega}\delta_{H(t,\omega)}(A_{\theta_t\omega})P(d\omega).	
	\end{equation}
\end{thm}
\noindent {\bf Proof.}\ \ It's obvious that P is the marginal measure of $\mu_s$ on $(\Omega,\cal F)$. Thus $\mu_s\in{\cal P}_P(\Omega\times X)$.
$\forall t\in\mathbb {R}, A\in{\cal F}\times{\cal B}(X)$,
\begin{eqnarray}\label{55}
({\Theta}^{-1}(t)(A))_\omega
&=&\{x:(\theta_t\omega,\Phi(t,\omega)x)\in A\}\nonumber\\
&=&\{x:\Phi(t,\omega)x\in A_{\theta_t}\omega\}\nonumber\\
&=&\Phi^{-1}(t,\omega)A_{\theta_t}\omega
\end{eqnarray}
So by (\ref{52}) and the definition of strong random almost periodic solution, we have
\begin{eqnarray}\label{56}
\Theta_t\mu_s(A)
&=&\mu_s({\Theta}^{-1}_t(A))\nonumber\\
&=&\int_{\Omega}\delta_{H(s,\omega)}({\Theta}^{-1}_t(A))_{\theta_s\omega}P(d\omega)\nonumber\\
&=&\int_{\Omega}\delta_{H(s,\omega)}(\Phi^{-1}(t,\theta_s\omega)A_{\theta_t\theta_s\omega})P(d\omega)\nonumber\\
&=&\int_{\Omega}\delta_{\Phi(t,\theta_s\omega)H(s,\omega)}(A_{\theta_{t+s}\omega})P(d\omega)\nonumber\\
&=&\int_{\Omega}\delta_{H(t+s,\omega)}(A_{\theta_{t+s}\omega})P(d\omega)\nonumber\\
&=&\mu_{t+s}(A).
\end{eqnarray}
By the definition of random almost periodic solution, for any $\varepsilon>0$, there exists a number $l(\varepsilon)$ defining a relatively dense set of numbers $\{\tau_k\}$ shch that for almost all $\omega\in\Omega$, $d(H(s+\tau_k, \omega),H(s,\theta(\tau_k)\omega))<\varepsilon, \ \forall t,s\in\mathbb{R}$ and  for any $s\in \mathbb{R}$, $\tau_k\in\{\tau_k\}$ by the probability preserving property of $\theta$, we have
\begin{eqnarray}\label{57}
\rho(\mu_{s+\tau_k},\mu_s)
&=&sup|\int_{\Omega}f(\theta_s\omega,H(s,\omega))P(d\omega)-\int_{\Omega}f(\theta_{s+\tau_k}\omega,H(s+\tau_k,\omega))P(d\omega)|\nonumber\\
&=&sup|\int_{\Omega}f(\theta_{s+\tau_k}\omega,H(s,\theta_{\tau_k}\omega))P(d\omega)-\int_{\Omega}f(\theta_{s+\tau_k}\omega,H(s+\tau_k,\omega))P(d\omega)|\nonumber\\
&\leq&sup\int_{\Omega}||f(\theta_{s+\tau_k}\omega,H(s,\theta_{\tau_k}\omega))-f(\theta_{s+\tau_k}\omega,H(s+\tau_k,\omega))||P(d\omega)\nonumber\\
&\leq&C_1d(H(s,\theta_{\tau_k}\omega),H(s+\tau_k,\omega))\nonumber\\
&=&C_1\varepsilon.
\end{eqnarray}
So by the definition of 5.1, $\mu_t, t\in\mathbb{R}$ is an almost periodic probability measure.
\hfill\fbox\\

\begin{rem}
	Let $\Phi$ be a measurable random dynamical system over $\theta$. Suppose $H$ is a random almost periodic solution of $\Phi$, it is easy to see that the factorization of $\mu_t$ defined in Theorem 5.2 is
	$$(\mu_t)_\omega=\delta_{H(t,\theta_{-t}\omega)},\ \forall t\in\mathbb{R}.$$
\end{rem}

In order to obtain the almost periodic probability measureon the phase space, we consider a Markovian cocycle random dynamical system $\Phi:\mathbb{R}^+\times\Omega\times X\rightarrow X$ on a filtered dynamical system $(\Omega,{\cal F},P,(\theta_t)_{t\in\mathbb{R}},({\cal F}^s_t)_{t\leq s})$. Denote the transition probability of Markovian process $\Phi(t,\omega)x$ on the Polish space $(X,{\cal B}(X))$ by (cf.\cite{A,GJ})
$$P(t,x,B)=P(\{\omega:\Phi(t,\omega)x\in B\}),\ \forall t\in\mathbb{R}^+, B\in{\cal B}(X),$$
and for any $t\geq 0$, $\lambda\in{\cal P}(X)$ we set
$$P_t^*\lambda(B)=\int_X P(t,x,B)\lambda(dx), \ t\geq 0,\ B\in{\cal B}(X).$$

\begin{defi}
	Given a Markovian random dynamical system $\Phi:\mathbb{R}^+\times\Omega\times X\rightarrow X$ over $\theta$, a measure $\lambda:\mathbb{R}\rightarrow{\cal P}(X)$ is called an almost periodic probability measure on the phase space $(X, {\cal B}(X))$ if for any $\varepsilon>0$, there exists a number $l(\varepsilon)$ defining a relatively dense set of numbers $\{\tau_k\}$ such that for almost all $\omega\in\Omega$,
	\begin{equation}\label{58}
	P_t^*\lambda_s=\lambda_{t+s}, \ \rho_1(\lambda_{{\tau_k}+s},\lambda_s)<\varepsilon, \ \forall t\in\mathbb{R}^+, s\in\mathbb{R}.
	\end{equation}
\end{defi}

\begin{thm}
	Let $\Phi:\mathbb{R}^+\times\Omega\times X\rightarrow X$ be a Markovian random dynamical system over $\theta$, if it has an adapted random  almost periodic solution $H:\mathbb{R}\times\Omega\rightarrow X$ with relatively dense set $\{\tau_k\}$. Then it has an almost periodic measure $\lambda:\mathbb{R}\rightarrow{\cal P}(X)$ on the phase space $(X, {\cal B}(X))$. For any $g\in C_{BL}(X)$,
	\begin{equation}\label{59}
	\lambda_t(g)=\int_{\Omega}g(H(t,\theta_{-t}\omega))P(d\omega)
	\end{equation}	
	and for any $A\in{\cal B}(X)$,	
	\begin{equation}\label{60}
	\lambda_t(A)=P(\{\omega:H(t,\theta_{-t} \omega)\in A\}).	
	\end{equation}
\end{thm}
\noindent {\bf Proof.}\ \ Firstly, for any $A\in{\cal B}(X)$, $t\in\mathbb{R}^+$, and $s\in\mathbb{R}$, by measure preserving property of $\theta$ and independency of $\Phi(t,\theta_s\omega)$ and ${\cal F}^s_{-\infty}$,
\begin{eqnarray}\label{61}
\lambda_{t+s}(A)
&=&P(\{\omega:H(t+s,\omega)\in A\})\nonumber\\
&=&P(\{\omega:\Phi(t,\theta_s\omega)H(s,\omega)\in A \})\nonumber\\
&=&\int_X P(t,x,A)P(\{\omega:H(s,\omega)\in dx\})\nonumber\\
&=&\int_X P(t,x,A)\lambda_s(dx)\nonumber\\
&=&P_t^*\lambda_s(A).
\end{eqnarray}
Secondly, for any $s\in\mathbb{R}$, $\tau_k\in\{\tau_k\}$ and almost all $\omega\in\Omega$, by (\ref{50}) and the definition of random almost periodic solution we have:
\begin{eqnarray}\label{62}
\rho_1(\lambda_{s+\tau_k},\lambda_s)
&=&sup|\int_{\Omega}g(H(s,\theta_{-s}\omega))P(d\omega)-\int_{\Omega}g(H(s+\tau_k,\theta_{-s-\tau_k}\omega))P(d\omega)|\nonumber\\
&\leq&sup\int_{\Omega}||g(H(s,\theta_{-s}\omega))-g(H(s+\tau_k,\theta_{-s-\tau_k}\omega))||P(d\omega)\nonumber\\
&\leq&C_2d(H(s,\theta_{-s}\omega),H(s+\tau_k,\theta_{-s-\tau_k}\omega))\nonumber\\
&=&C_2\varepsilon.
\end{eqnarray}
So by the definition of 5.4, $\lambda_t, t\in\mathbb{R}$ is an almost periodic probability measure on the phase space $X$.
\hfill\fbox\\

{ \noindent {\bf\large Acknowledgements}\ \
	The authors would like to thank Wei Sun for helpful discussions about this paper. Z.H. Zheng acknowledges financial supports of the NSF of China (No. 11671382), CAS Key Project of Frontier Sciences (No. QYZDJ-SSW-JSC003), the Key Lab. of Random Complex Structures and Data Sciences CAS and National Center for Mathematics and Interdisciplinary Sciences CAS.

\end{document}